\documentclass[letterpaper, 10 pt, conference]{ieeeconf}

\IEEEoverridecommandlockouts

\let\labelindent\relax

\usepackage{cite}
\usepackage[pdftex]{graphicx}
\usepackage[cmex10]{amsmath}
\usepackage{amsthm}
\usepackage{amssymb}
\usepackage{enumitem}
\usepackage{savesym}
\savesymbol{AND}
\usepackage{algorithmic}
\usepackage{algorithm}
\usepackage{tikz}
\usepackage{caption}
\usepackage{subcaption}
\usepackage{array}
\usepackage{booktabs}

\usetikzlibrary{decorations.pathreplacing,shapes,arrows,calc}

\tikzstyle{block} = [draw, fill=white!20, rectangle, 
    minimum height=1cm, minimum width=1cm,rounded corners=5pt]
\tikzstyle{sum} = [draw, fill=white!20, circle, node distance=1cm]
\tikzstyle{input} = [coordinate]
\tikzstyle{output} = [coordinate]
\tikzstyle{tmp} = [coordinate]
\tikzstyle{pinstyle} = [pin edge={<-,line width=0.3mm,black}]

\theoremstyle{plain}
\newtheorem{thm}{Theorem}

\theoremstyle{definition}
\newtheorem{defn}{Definition}

\newtheorem{assum}{Assumption}

\newcommand{\algorithmiconline}{\textbf{Online: }}
\newcommand{\ONLINE}{\STATE \algorithmiconline}
\newcommand{\algorithmicoffline}{\textbf{Offline: }}
\newcommand{\OFFLINE}{\STATE \algorithmicoffline}

\begin{document}

\title{\LARGE \bf Tractable Dual Optimal Stochastic Model Predictive Control:\\
An Example in Healthcare}

\author{
  Martin A. Sehr \& Robert R. Bitmead
  \thanks{The authors are with the Department
    of Mechanical and Aerospace Engineering, UC San Diego, La
    Jolla, CA 92093, USA. 
\newline{\tt\small \{msehr,rbitmead\}@ucsd.edu}}%
}

\maketitle

\begin{abstract}
Output-Feedback Stochastic Model Predictive Control based on Stochastic Optimal Control for nonlinear systems is computationally intractable because of the need to solve a Finite Horizon Stochastic Optimal Control Problem. However, solving this problem leads to an optimal probing nature of the resulting control law, called dual control, which trades off benefits of exploration and exploitation. In practice, intractability of Stochastic Model Predictive Control is typically overcome by replacement of the underlying Stochastic Optimal Control problem by more amenable approximate surrogate problems, which however come at a loss of the optimal probing nature of the control signals. While probing can be superimposed in some approaches, this is done sub-optimally. In this paper, we examine approximation of the system dynamics by a Partially Observable Markov Decision Process with its own Finite Horizon Stochastic Optimal Control Problem, which can be solved for an optimal control policy, implemented in receding horizon fashion. This procedure enables maintaining probing in the control actions. We further discuss a numerical example in healthcare decision making, highlighting the duality in stochastic optimal receding horizon control.
\end{abstract}
\section*{Glossary}
\begin{description}
\item[FHSOCP] finite-horizon stochastic optimal control problem,
\item[HMM] hidden Markov model,
\item[i.i.d.] independent and identically distributed,
\item[POMDP] partially observable Markov decision process,
\item[SDPE] stochastic dynamic programming equation,
\item[SMPC] stochastic model predictive control,
\item[UAV] unmanned autonomous/aerial vehicle.
\end{description}

\section{Introduction}
\label{intro}
Stochastic Optimal Control on the infinite time horizon is a computationally intractable problem. Unfortunately, the same holds even after replacing the infinite-horizon control problem by receding horizon implementation of a finite-horizon stochastic optimal control law. While a number of approaches have been devoted to Stochastic Model Predictive Control (SMPC), including~\cite{yan2005incorporating,sui2008robust,mayne2009robust,IFAC2017}, none of these involve solution of a Finite Horizon Stochastic Optimal Control Problem (FHSOCP), which is required for optimal \emph{probing} in the resulting control inputs but computationally intractable in practice. Assuming an available solution of the FHSOCP in principle leads to a \emph{dual optimal} SMPC law, which enjoys a number of desirable properties discussed in~\cite{sehr2016stochastic}. These properties include recursive feasibility, stochastic stability and, in particular, bounds relating infinite-horizon performance of the SMPC control with computed finite-horizon performance and optimal infinite-horizon performance. This result, which requires explicit solution of the FHSOCP, encourages a tractable version of dual optimal SMPC.

In this paper, we discuss a version of SMPC motivated by the structure of dual optimal SMPC as in~\cite{sehr2016stochastic}. In particular, we suggest approximation of the system dynamics by a Partially Observable Markov Decision Process (POMDP), the finite-horizon solution of which is tractable for small to moderate problem instances (see e.g.~\cite{smallwood1973optimal,kaelbling1998planning}). The main benefit of this approach is that duality of the resulting control inputs is naturally maintained by the approximation in an optimal sense. Moreover, given that POMDPs are a subclass of general nonlinear systems, the results in~\cite{sehr2016stochastic} hold in modified form for the resulting approximate dynamics. These two observations lead us to propose SMPC on the approximate POMDP with explicit solution of the resulting approximate FHSOCP. Receding horizon control in POMDPs has been discussed for instance with respect to UAV guidance~\cite{miller2009pomdp}, sensor scheduling in~\cite{sunberg2013information}, and robotic manipulation in~\cite{pajarinen2015robotic}. We discuss in particular an example in medical decision making, involving decisions regarding appointment scheduling, the use of costly diagnostic tests, and patient treatment. This example highlights in particular the probing nature of the dual optimal SMPC on a POMDP. 

The paper evolves as follows. After briefly introducing SMPC as suggested in~\cite{sehr2016stochastic} in Section~\ref{sec:smpc}, we transition to receding horizon control of POMPDs in Section~\ref{sec:pomdp}. We transition to the use of POMDPs in healthcare in Section~\ref{sec:eg}, which includes our numerical example highlighting duality in receding horizon implementations of optimal POMDP solutions. The paper closes with some concluding remarks in Section~\ref{sec:conclusions}.

\section{Stochastic Model Predictive Control}
\label{sec:smpc}
\subsection{Stochastic Optimal Control in a Nutshell}
\begin{figure}[tb]
\begin{center}
\begin{tikzpicture}[auto, node distance=1cm,>=latex']
   \node [line width=0.3mm,block, pin={[pinstyle]above:$x_0,v_t,w_t$},
           node distance=3cm] (sys) {System~(\ref{eq:state}-\ref{eq:output})}; 
   \node [output, right of=sys,node distance=2cm] (output) {};
   \node [output, right of=output,node distance=1cm] (out2) {};
   \node [line width=0.3mm,block, below of=sys,node distance=2.6cm, pin={[pinstyle]above:$\pi_{0\mid -1}= \operatorname{pdf}(x_0)$}] (filter)
{Filter~(\ref{eq:BF_rec}-\ref{eq:BF_pred})};
   \node [line width=0.3mm,block, left of=filter, node distance=3cm]
   (ctr) {Control~\eqref{eq:policies}};
   \node [input, left of=ctr,node distance=2cm] (input) {};

   \draw [line width=0.3mm,->] (input) |- (sys);
   \draw [line width=0.3mm,-] (sys) -- (output);
   \draw [line width=0.3mm,->] (filter) -- node [above,name=pi]{$\pi_t$}(ctr);
   \draw [line width=0.3mm,->] (output) |- node [pos=0.75,swap]{$y_t$} (filter);
   \draw [line width=0.3mm,-] (ctr) -- node [above]{$u_t$} (input);

   \node [line width=0.3mm,block, below of=pi,node distance=1.5cm] (delay) {$z^{-1}$};
   
   \draw [line width=0.3mm,->] (delay) -| node [pos=0.25]{$u_{t-1}$} (filter);
   \draw [line width=0.3mm,->] (input) |- (delay);

\node (rect) at (-1.5,-3.4) [line width=0.4mm,draw,dotted,minimum width=7.5cm,minimum height=3.1cm] {};
\node[anchor=south west] at (rect.south west) {\emph{Output-Feedback Controller}};
\end{tikzpicture}
\caption{Closed-loop system architecture for stochastic optimal output-feedback control. The optimal control law is a \textit{separated} feedback policy, with information state update via Bayesian Filter and an information state feedback law, mapping information states $\pi_t$ to feasible control inputs $u_t$.}
 \label{fig:architecture}
\end{center}
\end{figure}
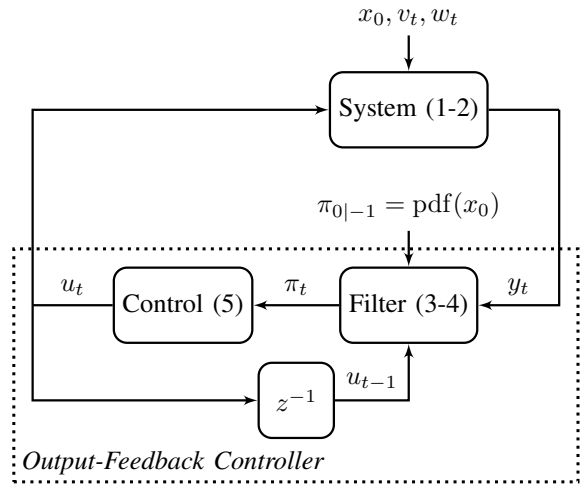

We consider receding horizon output-feedback control for nonlinear stochastic systems of the form
\begin{align}
x_{t+1}&=f(x_t,u_t,w_t),\quad x_0\in\mathbb{R}^n,\label{eq:state}\\
y_t&=h(x_t,v_t), \label{eq:output}
\end{align}
where $x_k\in\mathbb{R}^n$, $u_k\in\mathbb{R}^m$, $y_k\in\mathbb{R}^o$, starting from known initial state density $\pi_{0|-1} = \operatorname{pdf}(x_0)$. We denote the data available at time $t$ by 
\begin{align*}
\mathbf{\zeta}^t&\triangleq\{y_0,u_0,y_1,u_1,\dots,u_{t-1},y_t\},&
\mathbf{\zeta}^0&\triangleq\{y_0\} .
\end{align*}
We further impose the following standing assumption on the random variables and control inputs.
\begin{assum}\label{assm:sys}
The signals in~(\ref{eq:state}-\ref{eq:output}) satisfy:
\begin{enumerate}[label=\arabic*.]
\item $w_t$ and $v_t$ are i.i.d. sequences with known densities.
\item $x_0, w_t, v_l$ are mutually independent for all $t,l\geq 0$.
\item The control input $u_t$ at time instant $t\geq 0$ is a function of the data $\mathbf{\zeta}^t$ and given initial state density $\pi_{0\mid -1}$.
\end{enumerate}
\end{assum}
The \textit{information state,} denoted $\pi_t$, is the conditional probability density function of state $x_t$ given data $\mathbf{\zeta}^t$,
\begin{align*}
\pi_{t}&\triangleq\operatorname{pdf}\left(x_{t}\mid \mathbf{\zeta}^t \right),
\end{align*}
and is propagated via the \emph{Bayesian Filter} (see e.g.~\cite{chen2003bayesian,simon2006optimal}):
\begin{align}
\pi_{t}&=
\frac{\operatorname{pdf}(y_{t}\mid x_{t})\,\pi_{t\mid t-1}}{\int \operatorname{pdf}(y_{t}\mid x_{t})\,\pi_{t\mid t-1}\,dx_{t}},\label{eq:BF_rec} \\
\pi_{t+1\mid t} &\triangleq
\int \operatorname{pdf}(x_{t+1} \mid x_{t},u_{t}) \,\pi_{t}\, dx_{t},\label{eq:BF_pred}
\end{align}
from initial density $\pi_{0\mid -1}$. As a result of the Markovian dynamics~(\ref{eq:state}-\ref{eq:output}), optimal control inputs must inherently be \textit{separated} feedback policies (e.g.~\cite{bertsekas1995dynamic,BKKUM1986}). That is, control input $u_{t}$ depends on the data $\mathbf{\zeta}^t$ and initial density $\pi_{0\mid -1}$ solely through the current information state, $\pi_{t}$. Optimality thus requires propagating $\pi_{t}$ and policies $g_t$, where
\begin{align}\label{eq:policies}
u_t = g_t(\pi_{t}),
\end{align}
leading to the closed-loop architecture displayed in Figure~\ref{fig:architecture}.
\begin{defn}
$\mathbb{E}_k[\,\cdot\,]$ and $\mathbb{P}_k[\,\cdot\,]$ are expected value and probability with respect to state $x_k$ -- with conditional density $\pi_k$ -- and i.i.d. random variables $\{(w_j,v_{j+1}):j\geq k\}$.
\end{defn}
Given control horizon $N\in\mathbb{N}$, stage cost $c:\mathbb{R}^{n}\times\mathbb{R}^{m}\to\mathbb{R}_+$, terminal cost $c_N:\mathbb{R}^{n}\to\mathbb{R}_+$, and discount factor $\alpha\in\mathbb{R}_+$, we further define the cost
\begin{align*}
J_N(\pi_{0},\mathbf{g}^{N-1})= 
\mathbb{E}_0\left[\sum_{j=0}^{N-1}{\alpha^j c(x_j,g_j(\pi_{j}))} + \alpha^{N}c_{N}(x_{N})\right],
\end{align*}
where $\mathbf{t}^{m} \triangleq \{t_0,t_{1},\ldots,t_{m}\}$. Extending this cost to the \emph{infinite horizon} typically requires $\alpha < 1$ and omitting the terminal cost term $c_N(\cdot)$, leading to
\begin{align*}
J_\infty(\pi_{0},\mathbf{g}^{\infty})&\triangleq
\mathbb{E}_0\left[\sum_{j=0}^\infty{\alpha^jc(x_j,g(\pi_{j}))}\right].
\end{align*}
Additionally introducing the constraints
\begin{align*}
\mathbb{P}_k\left[ x_k \in \mathcal{X}_k \right] &\geq 1 - \epsilon_k,&k&\in\mathbb{N}_1,\\
u_k = g_k(\pi_k) &\in \mathcal{U}_k, &k&\in\mathbb{N}_0,
\end{align*}
for $\epsilon_k \in[0,1)$, finite-horizon stochastic optimal feedback policies may be computed, in principle, by solving the Stochastic Dynamic Programming Equation (SDPE),
\begin{equation}
\begin{aligned}
\label{eq:DP1}
V_k(\pi_{k}) \triangleq 
\inf_{g_k(\cdot)}&\ \mathbb{E}_k \left[ c(x_k,g_k(\pi_k)) +\alpha V_{k+1}(\pi_{k+1})\right], \\
\text{s.t.}\,& \ \mathbb{P}_{k+1}\left[ x_{k+1} \in \mathbb{X}_{k+1} \right]\geq 1 - \epsilon_{k+1},\\
& \ g_k(\pi_k) \in \mathbb{U}_k,
\end{aligned}
\end{equation}
for $k = 0,\ldots,N-1$. The equation is solved backwards in time, from its terminal value
\begin{align}\label{eq:DP2}
V_N(\pi_{N}) &\triangleq \mathbb{E}_N \left[ c_N(x_{N})\right].
\end{align}

Solution of the SDPE is the primary source of intractability in Stochastic Optimal Control. The reason for this difficulty lies in the dependence of the future information states in~(\ref{eq:DP1}-\ref{eq:DP2}) on current and future control inputs via~(\ref{eq:BF_rec}-\ref{eq:BF_pred}). However, optimality via the SDPE leads to a control law of \emph{dual} nature~\cite{Feldbaum1965}. 
Notice that, while the Bayesian Filter~(\ref{eq:BF_rec}-\ref{eq:BF_pred}) can be approximated to arbitrary accuracy using a Particle Filter \cite{simon2006optimal}, the SDPE cannot be easily simplified without loss of optimal probing in the control inputs. While control laws generated without solution of the SDPE can be modified artificially to include certain excitation properties, as discussed for instance in~\cite{genceli1996new,marafioti2014persistently}, such approaches are suboptimal and do not enjoy the same theoretical guarantees.

\subsection{Dual Optimal SMPC Algorithm}
While computationally intractable, solution of the SDPE and the associated feedback architecture in Figure~\ref{fig:architecture} lead to the following dual optimal SMPC algorithm, the properties of which are discussed in \cite{sehr2016stochastic}.
\begin{algorithm}[H] 
\caption{Dual Stochastic Model Predictive Control}\label{RHSOC}                       
\begin{algorithmic}[1]
\OFFLINE
\STATE 
Solve~(\ref{eq:DP1}-\ref{eq:DP2}) for first optimal policy, $g_0^{\star}(\cdot)$
\ONLINE
\FOR{$t=0,1,2,\ldots$}
\STATE Measure $y_t$
\STATE Compute $\pi_t$
\STATE Apply first optimal control policy, $u_{t} = g_0^{\star}(\pi_{t})$
\STATE Compute $\pi_{t+1\mid t}$
\ENDFOR
\end{algorithmic}
\end{algorithm}
This algorithm differs from common practice in stochastic model predictive control in that it explicitly uses the information state $\pi_t$. Throughout the literature, these conditional densities are commonly replaced by state estimates in order to generate certainty-equivalent control laws. While this makes the problem more tractable, one no longer solves the underlying Stochastic Optimal Control problem, leading to a loss of duality in the resulting control inputs. The central divergence to such approaches lies in Step~2 of the algorithm, in which the SDPE is presumed solved offline for the optimal feedback policies $g_t^{\star}(\pi_t)$. This is an extraordinarily difficult proposition in many cases but captures the optimality, and hence duality, as a closed-loop feedback control law. The complexity of this step lies not only in computing a vector functional but also in the internal propagation of the information state within the SDPE.
\begin{defn}\label{def:opt}
Denote by $\mathbf{g}^{MPC}$ the SMPC implementation of policy $g_0^{\star}(\cdot)$ on the infinite horizon, i.e. 
\begin{align*}
\mathbf{g}^{MPC} \triangleq \{ g_0^{\star},g_0^{\star},g_0^{\star},\ldots\}.
\end{align*}
Similarly, $\mathbf{g}^{\star^{N-1}}$ and $\mathbf{g}^{\star^\infty}$ are the optimal sequences of policies in~(\ref{eq:DP1}-\ref{eq:DP2}) and the corresponding infinite-horizon problem.
\end{defn}
Subsuming assumptions akin to those used in deterministic MPC (e.g.~\cite{grune2008infinite}) to ensure stochastic recursive feasibility of the SMPC algorithm and constrain the growth of the underlying value function with increasing control horizon $N$, we re-state the following technical result from~\cite{sehr2016stochastic}.
\begin{thm}\label{thm:bounds}
Under suitable assumptions, SMPC with discount factor $\alpha\in[0,1)$ yields
\begin{multline*}
(1- \alpha\gamma)\, J_{\infty}(\pi_0,\mathbf{g}^{\star^\infty}) \leq 
(1- \alpha\gamma)\, J_{\infty}(\pi_0,\mathbf{g}^{MPC}) \\ \leq 
J_{N}(\pi_0,\mathbf{g}^{\star^{N-1}}) + \frac{\alpha}{1 - \alpha}\eta,
\end{multline*}
for given values of $\gamma\in[0,1]$ and $\eta\in\mathbb{R}_+$.
\end{thm} 
This result relates the following quantities: \textit{design cost}, $J_N(\pi_0,\mathbf{g}^{*^{N-1}}),$ which is known as part of the SMPC calculation, \textit{optimal cost} $J_\infty(\pi_0,\mathbf{g}^{*^\infty})$ which is unknown (otherwise we would use $\mathbf{g}^{*^\infty}$), and unknown infinite-horizon SMPC \textit{achieved cost} $J_\infty(\pi_0,\mathbf{g}^{MPC})$. The result, which requires duality and solution of the SDPE, is special in that SMPC approaches relying of approximation of the SDPE, as commonly found in the literature, do generally not yield statements regarding performance of the implemented control law on the infinite horizon.

\section{Approximate Dual Optimal SMPC via POMDPs}
We next approximate the nonlinear dynamics~(\ref{eq:state}-\ref{eq:output}) by a POMDP. The above algorithm and observations are then applicable in slightly modified form for the resultant approximate model. However, we can solve finite-horizon POMDP problems of reasonable dimensions explicitly, allowing us to maintain duality in the approximate version of the SMPC algorithm in Section~\ref{sec:smpc}. For approximation of the system dynamics and cost via POMDPs, see for example~\cite{sunberg2013information}, which discusses a related approach. In contrast with this reference, we suggest choosing the approximation so that the resulting SDPE can be solved explicitly. This results in maintaining duality in a stochastic optimal sense, although on the approximate POMDP dynamics.

\label{sec:pomdp}
\subsection{Partially Observable Markov Decisions Processes}
POMDPs are characterized by probabilistic dynamics on a finite state space $\mathbf{X} = \{1,\ldots,n\}$, finite action space $\mathbf{U} = \{1,\ldots,m\}$, and finite observation space $\mathbf{Y} = \{1,\ldots,o\}$. POMDP dynamics are defined by the conditional state transition and observation probabilities
\begin{align}\label{eq:pomdpx}
\mathbb{P}\left(x_{t+1} = j \mid x_t = i, u_t = a\right) &= p_{ij}^{a}, \\
\label{eq:pomdpy}
\mathbb{P}\left(y_{t+1} = \theta \mid x_{t+1} = j, u_t = a\right) &= r_{j\theta}^{a},
\end{align}
$t \geq 1$, $i,j\in\mathbf{X}$, $a\in\mathbf{U}$, $\theta\in\mathbf{Y}$. The state transition dynamics~\eqref{eq:pomdpx} correspond to a conventional Markov Decision Process (MDP). However, the control actions $u_t$ are to be chosen based on the initial state distribution $\pi_0 = \operatorname{pdf}(x_0)$ and the sequences of observations, $\{y_1,\ldots,y_t\}$, and controls $\{u_0,\ldots,u_{t-1}\}$, respectively. That is, we are choosing our control actions in a Hidden Markov Model (HMM~\cite{elliott2008hidden}) setup.

Given control action $u_t = a$ and measured output $y_{t+1} = \theta$, the Bayesian Filter recursion~(\ref{eq:BF_rec}-\ref{eq:BF_pred}) extends to the POMDP dynamics~(\ref{eq:pomdpx}-\ref{eq:pomdpy}) as
\begin{align*}
\pi_{t+1,j} = \frac{\sum_{i\in\mathbf{X}} 
\pi_{t,j} p_{ij}^a r_{j\theta}^a}{\sum_{i,j\in\mathbf{X}} \pi_{t,j} p_{ij}^a r_{j\theta}^a} ,
\end{align*}
where $\pi_{t,j}$ denotes the $j^{\text{th}}$ entry of the row vector $\pi_t$. We define the cost as in Section~\ref{sec:smpc}, with stage cost $c(x_t,u_t) = c_{i}^{a}$ if $x_t = i\in\mathbf{X}$ and $u_t = a\in\mathbf{U}$, or $e_ic(a)$ in vectorized form. The terminal cost $c_N(x_t)$ is defined similarly. Omitting constraints for brevity, optimal control decisions may then be computed via SDPE
\begin{multline*}
V_k(\pi_k) = \\
\min_{u_k\in\mathbf{U}} \left\{ \pi_{k} c(u_k) + 
\alpha\sum_{\theta \in\mathbf{Y}} \mathbb{P}\left(\theta \mid \pi_k,\,u_k \right) 
V_{k+1}(\pi_{k+1}) \right\}, 
\end{multline*}
from terminal value function
\begin{align*}
V_N(\pi_N) = \pi_N c_N .
\end{align*}

\section{Application in Healthcare}
\label{sec:eg}
\subsection{POMDPs in Healthcare Decision Making}
A framework for the use of POMDPs in healthcare decision-making applications is discussed in detail in~\cite{bennett2013artificial}. Particular examples of POMDPs capturing disease stages as well as appropriate tests and treatment decisions for individual patients include~\cite{hauskrecht2000planning,ayer2012or}. In~\cite{ayer2012or}, for instance, a POMDP is used for mammography screening in individual patients. States in the underlying Markov chain model are
\begin{enumerate}
\item No cancer
\item In situ cancer
\item Invasive cancer
\item In situ post-cancer
\item Invasive post-cancer
\item Death
\end{enumerate}
Available decisions are either to perform a mammogram, with possible follow-up biopsy, or to wait, which may lead to self-detection. The authors proceed with solution of the POMDP for optimal screening strategies depending on patient age and risk for in situ or invasive cancer. This highlights a successful application of POMDPs in healthcare.

\subsection{An Illustrative POMDP in Healthcare}
The remainder of this paper discusses a particular numerical example of decisions on treatment and diagnosis in healthcare, displaying specifically the use of dual control in SMPC applied to a POMDP. Consider a patient treated for a specific disease which can be managed but not cured. For simplicity, we assume that the patient does not die under treatment. While this transition would have to be added in practice, it results in a time-varying model, which we avoid in order to keep the following discussion compact. 

\begin{figure}[tb]
  \centering
  \includegraphics[width=\columnwidth]{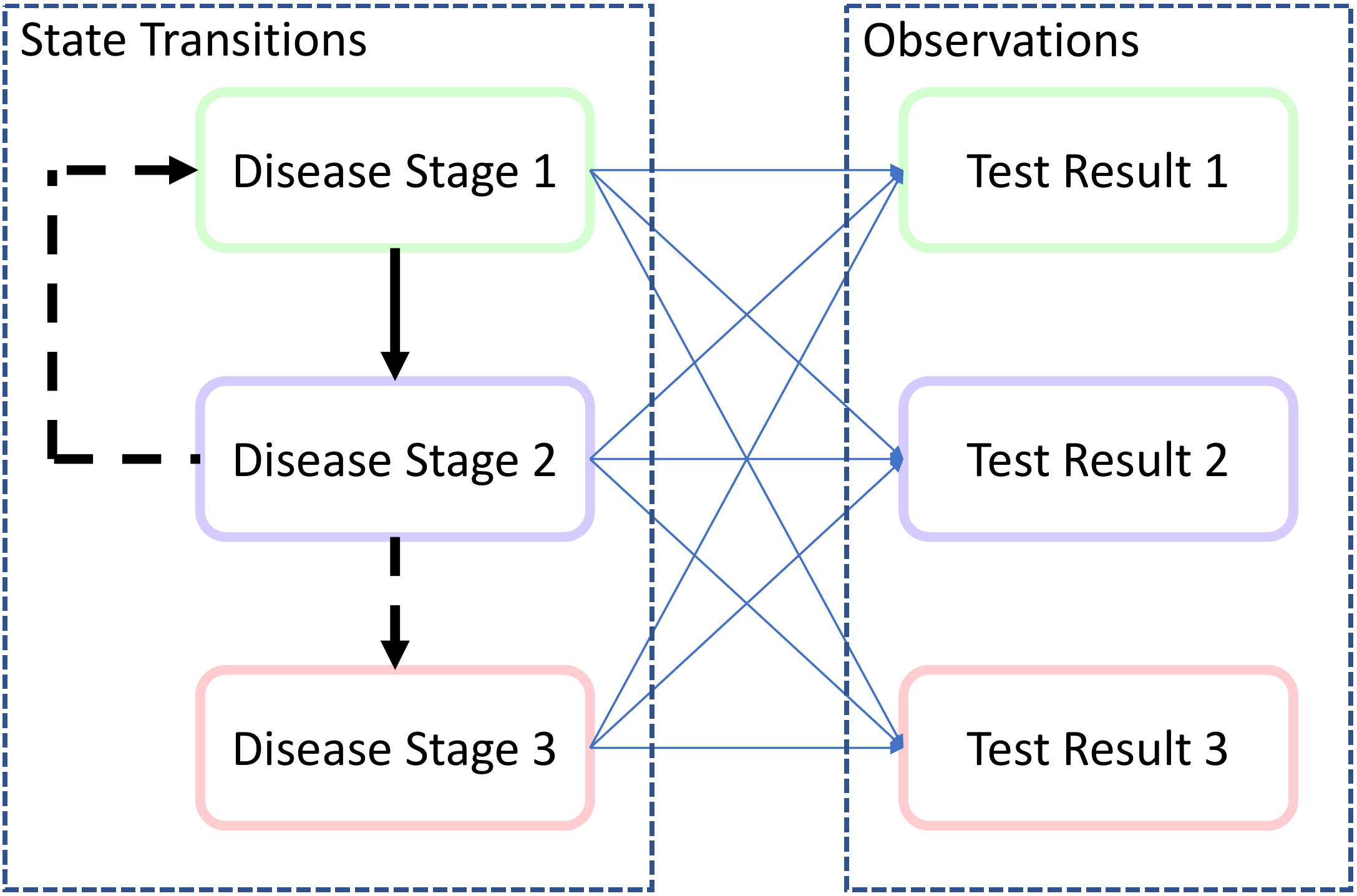}
  \caption{Feasible state transitions and possible test results in healthcare example. Solid arrows for feasible state transitions and observations. Dashed arrows for transitions conditional on treatment and diagnosis decisions.}
  \label{fig:transitions}
\end{figure}

The disease encompasses three stages with severity increasing from Stage 1 through Stage 2 to Stage 3, transitions between which are governed by a Markov chain with transition probability matrix
\begin{align*}
P = \begin{bmatrix}
0.8&0.2&0.0\\0.0&0.9&0.1\\0.0&0.0&1.0
\end{bmatrix},
\end{align*}
where $P$ is the matrix with values $p_{ij}$ at row $i$ and column $j$. All transition and observation probability matrices below are defined similarly. Once our patient enters Stage 3, Stages 1 and 2 are inaccessible for all future times. However, Stage 3 can only be entered through Stage 2, a transition from which to Stage 1 is possible only under costly treatment. The same treatment inhibits transitions from Stage 2 to Stage 3. We have access to the patient state only through tests, which will result in one of three possible values, each of which is representative of one of the three disease stages. However, these tests are imperfect, with non-zero probability of returning an incorrect disease stage. All possible state transitions and observations are illustrated in Figure~\ref{fig:transitions}.

At each point in time, the current information state $\pi_t$ is available to make one of four possible decisions:
\begin{enumerate}
\item Skip next appointment slot
\item Schedule new appointment
\item Order rapid diagnostic test
\item Apply available treatment
\end{enumerate}
Skipping an appointment slot results in the patient progressing through the Markov chain describing the transition probabilities of the disease without medical intervention, without new information being available after the current decision epoch. Scheduling an appointment does not alter the patient transition probabilities but provides a low-quality assessment of the current disease stage, which is used to refine the next information state. The third option, ordering a rapid diagnostic test, allows for a high-quality assessment of the patient's state, leading to a more reliable refinement of the next information state than possible when choosing the previous decision option. The results from this diagnostic test are considered available sufficiently fast so that the patient state remains unchanged under this decision. The remaining option entails medical intervention, allowing transition from Stage 2 to Stage 1 while preventing transition from Stage 2 to Stage 3. Transition probabilities $P(a)$, observation probabilities $R(a)$, and stage cost vectors $c(a)$ for each decision are summarized in Table~\ref{tab:numbers}. Additionally, we impose the terminal cost
\begin{align*}
c_N = \begin{bmatrix} 0 & 4 & 30 \end{bmatrix}^T.
\end{align*}
While we choose the SMPC discount factor as $\alpha = 1$ for simplicity, one may also extend this discussion and choose a discount factor $\alpha\in(0,1)$ to exploit Theorem~\ref{thm:bounds} explicitly. Given that we do not consider constraints in this example, such that recursive feasibility holds by default, this can be enabled via proper choice of the terminal cost term, $c_N$.

\begin{table*}[tb]
  \caption{Problem data for healthcare decision making example.}
  \centering
  \begin{tabular}{l|ccc}	
    Decision $a$ & Transition Probabilities $P(a)$ & Observation Probabilities $R(a)$ & Cost $c(a)$ \\
    \midrule
\hspace{0.2cm}1: Skip next appointment slot & 
$\begin{bmatrix} 0.80&0.20&0.00\\0.00&0.90&0.10\\0.00&0.00&1.00 \end{bmatrix}$ & 
$\begin{bmatrix} 1/3&1/3&1/3\\1/3&1/3&1/3\\1/3&1/3&1/3 \end{bmatrix}$ & 
$\begin{bmatrix} 0\\5\\5 \end{bmatrix}$ \\[0.5cm]
\hspace{0.2cm}2: Schedule new appointment & 
$\begin{bmatrix} 0.80&0.20&0.00\\0.00&0.90&0.10\\0.00&0.00&1.00 \end{bmatrix}$ & 
$\begin{bmatrix} 0.40&0.30&0.30\\0.30&0.40&0.30\\0.30&0.30&0.40 \end{bmatrix}$ & 
$\begin{bmatrix} 1\\1\\1 \end{bmatrix}$ \\[0.5cm]
\hspace{0.2cm}3: Order rapid diagnostic test & 
$\begin{bmatrix} 1.00&0.00&0.00\\0.00&1.00&0.00\\0.00&0.00&1.00 \end{bmatrix}$ & 
$\begin{bmatrix} 0.90&0.05&0.05\\0.05&0.90&0.05\\0.05&0.05&0.90 \end{bmatrix}$ & 
$\begin{bmatrix} 4\\3\\4 \end{bmatrix}$ \\[0.5cm]
\hspace{0.2cm}4: Apply available treatment & 
$\begin{bmatrix} 0.80&0.20&0.00\\0.75&0.25&0.00\\0.00&0.00&1.00 \end{bmatrix}$ & 
$\begin{bmatrix} 0.40&0.30&0.30\\0.30&0.40&0.30\\0.30&0.30&0.40 \end{bmatrix}$ & 
$\begin{bmatrix} 4\\2\\4 \end{bmatrix}$
  \end{tabular}
  \label{tab:numbers}
\end{table*}

\subsection{Rationale for Duality}
Intuitively, we expect an efficient policy for this problem to attempt avoiding transitions to Stage 3 while managing the ressources required to schedule appointments, order tests, or apply medical intervention. This may, in principle, be achieved by a policy akin to the following structure:
\begin{enumerate}
\item Skip appointments when Stages 2 and 3 are unlikely.
\item Schedule appointments when Stages 2 and 3 are likely but the probability for Stage 2 is below some threshold.
\item Order diagnostic test if the probability of Stage 2 lies in a specific range.
\item Proceed with medical intervention if the probability of Stage 2 is high.
\end{enumerate}
While the optimal policy may be somewhat more intricate, this simple decision structure could be acceptable in practice. However, even this simple structure includes duality in the decisions, demonstrated by including the diagnostic test even though it does not alter the patient state. That is, this decision improves the quality of available information at a cost, also called \emph{exploration}. This improvement in the available information allows us to apply medical intervention at appropriate times, which is called \emph{exploitation}. 

\subsection{Computational Results}
\begin{figure*}[tb]
  \centering
  \includegraphics[width=160mm]{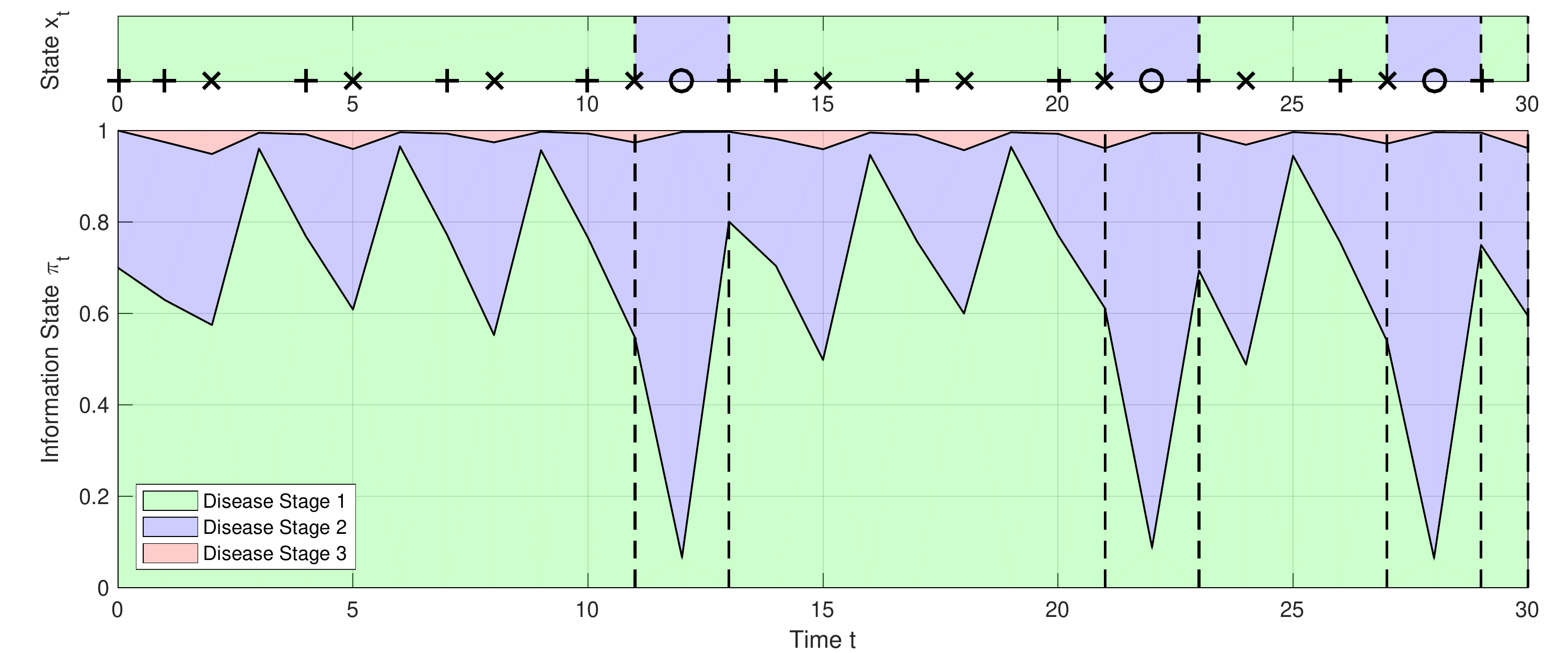}
  \caption{Simulation results for SMPC with horizon $N=6$ and discount factor $\alpha = 1$. Top plot displays patient state and transitions, with optimal SMPC decisions based on current information state: appointment (pluses); diagnosis (crosses); treatment (circles). Bottom plot shows information state evolution. Dashed vertical lines mark time instances of state transitions. }
  \label{fig:simN6}
\end{figure*}

The trade-off between these two principal decision categories is precisely what is encompassed by duality, which we can include in an optimal sense by solving the SDPE and applying the resulting initial policy in receding horizon fashion. This is demonstrated in Figure~\ref{fig:simN6}, which shows simulation results for SMPC with control horizon $N = 6$. As anticipated, the stochastic optimal receding horizon policy shows a structure not drastically different from the decision structure motivated above. In particular, diagnostic tests are used effectively to decide on medical intervention.

\begin{figure*}[tb]
  \centering
  \includegraphics[width=160mm]{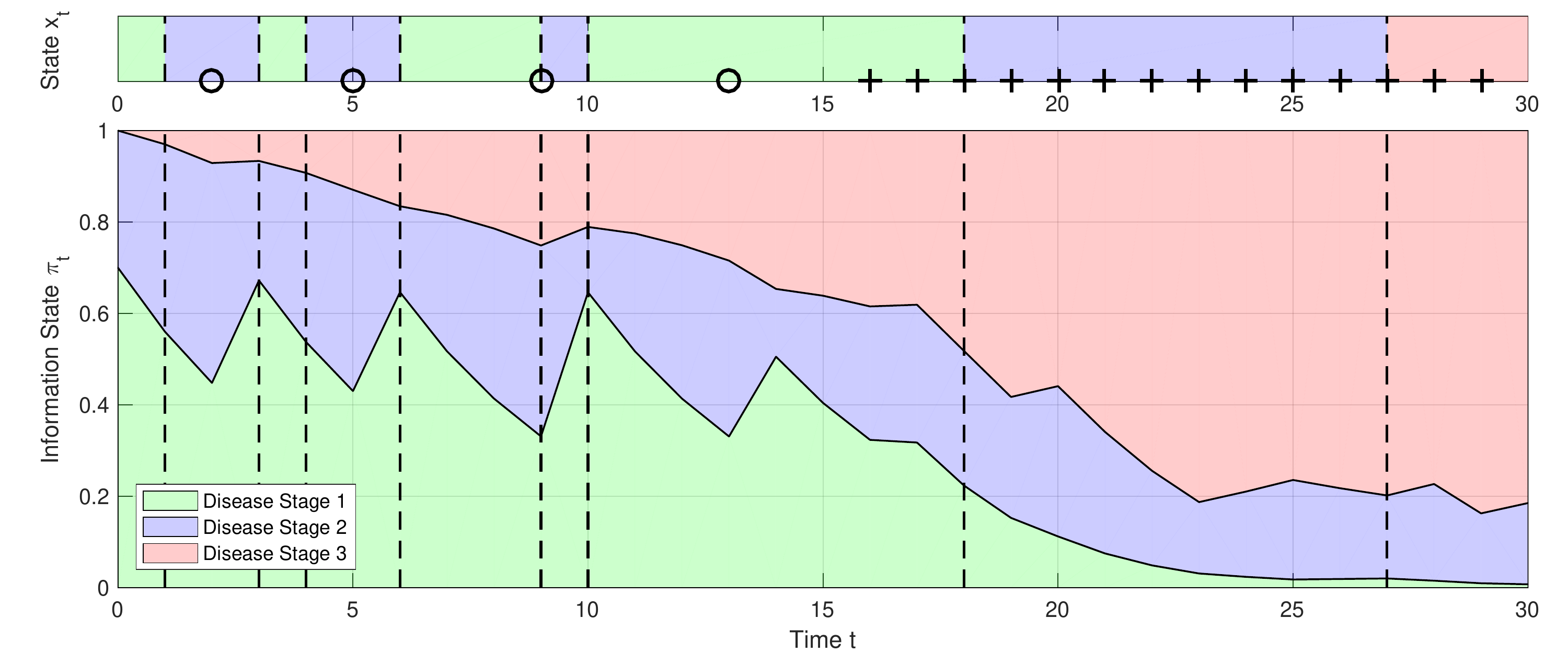}
  \caption{Simulation results for certainty-equivalent SMPC with horizon $N=6$ and discount factor $\alpha = 1$. Top plot displays patient state and transitions, with optimal SMPC decisions based on most likely patient state, extracted from current information state: appointment (pluses); diagnosis (crosses); treatment (circles). Bottom plot shows information state evolution. Dashed vertical lines mark time instances of state transitions.}
  \label{fig:simN6CE}
\end{figure*}

In contrast, Figure~\ref{fig:simN6CE} shows a similar simulation for certainty-equivalent SMPC with control horizon $N = 6$. In this control, the information state $\pi_t$ is replaced by the most likely patient state $\hat{x}_t$ at each time, a simplification often used to ease the computational burden associated with solving the underlying control problem and propagating the information state $\pi_t$ over time. However, this simplification comes at a loss of duality in the control decisions. As shown in Figure~\ref{fig:simN6CE}, the certainty equivalent control does not choose diagnosis to refine the information state and subsequently make a more informed treatment decision. Given that the control is only based on one of three possible estimate patient states, we can give the entire certainty equivalent optimal policy as:
\begin{enumerate}
\item Skip appointments when Stage 1 is most likely.
\item Apply treatment when Stage 2 is most likely.
\item Schedule appointments when Stage 3 is most likely.
\end{enumerate}
Clearly, this policy does not make effective use of the ressources and decision options available to us, caused by its loss of duality in the control signal. In fact, not classifying the patient state correctly towards the end of the simulation leads to a preventable transition to Stage 3 in this simulation. The failure of this seemingly similar certainty equivalent SMPC control strategy highlights the need for duality in the control signals, as is accommodated by dual optimal SMPC as motivated in this paper.

\section{Conclusions}
\label{sec:conclusions}
We discussed approximation of intractable nonlinear output-feedback SMPC problems by POMDPs, which for reasonable problem dimensions can be solved explicitly. The benefit of this approach as opposed to common practice in SMPC is that duality of the control actions natural to Stochastic Optimal Control is maintained optimally. We demonstrated the use of SMPC and in particular the probing nature of the resulting control inputs by means of a particular numerical POMDP in healthcare decision making. The example displays in particular how an expensive diagnostic test is chosen to refine the current information state without influencing the patient state. Without duality in the control input, this diagnostic test is not used, resulting subsequently in poorly informed treatment decisions due to lack of refinement in the information state.

\bibliographystyle{ieeetr}
\bibliography{CCTA2017}

\end{document}